\documentclass[11pt,reqno,oneside]{amsart}
\usepackage{amssymb}
\usepackage{graphics}
\usepackage[usenames]{color}
\usepackage[hidelinks]{hyperref}
\usepackage[utf8]{inputenc}

\newtheorem{lemma}{Lemma}[section]

\newtheorem{cor}{Corollary}[section]
\newtheorem{rem}{Remark}[section]

\newtheorem{cex}{Counterexample}[section]

\newcommand{\bit}{\begin{itemize}}
\newcommand{\eit}{\end{itemize}}
\newcommand{\bq}{\begin{equation}}
\newcommand{\eq}{\end{equation}}

\newcommand{\ds}{\displaystyle}

\newcommand{\bx}{\bf x}
\newcommand{\bv}{\bf v}

\newcommand{\R}{\mathbb R}

\makeatletter
\g@addto@macro{\endabstract}{\@setabstract}
\newcommand{\authorfootnotes}{\renewcommand\thefootnote{\@fnsymbol\c@footnote}}%
\makeatother

\title[Sharp geometric requirements in the error estimate]{Sharp geometric requirements in the Wachspress interpolation error estimate}



\author[G. Monz\'on]{Gabriel Monz\'on}
\address{Instituto de Ciencias\\ Universidad Nacional de General Sarmiento\\
	J. M. Guti\'errez 1150 \\ 
	(1613) Los Polvorines\\
	Buenos Aires\\ Argentina.}
\email{gmonzon@ungs.edu.ar}

\begin{document}
	
\maketitle

\begin{abstract}
	Geometric conditions on general polygons are given in \cite{GRB} in order to guarantee the error estimate for interpolants built from generalized barycentric coordinates, and the question about identifying sharp geometric restrictions in this setting is proposed. In this work, we address the question when the construction is made by using Wachspress coordinates. We basically show that the imposed conditions: {\it bounded aspect ratio property $(barp)$}, {\it maximum angle condition $(MAC)$} and {\it minimum edge length property $(melp)$} are actually equivalent to $[MAC,melp]$, and if any of these conditions is not satisfied, then there is no guarantee that the error estimate is valid. In this sense, $MAC$ and $melp$ can be regarded as sharp geometric requirements in the Wachspress interpolation error estimate.
\end{abstract}

\section{Introduction}

Many and different conditions on the geometry of finite elements were required in order to guarantee optimal convergence in the interpolation error estimate. Some of them deal with interior angles like the {\it maximum angle condition} (maximum interior angle bounded away from $\pi$) and the {\it minimum angle condition} (minimum interior angle bounded away from $0$), but others deal with some lengths of the element like the {\it minimum edge length property} (the diameter of the element is comparable to the length of the segment determined by any two vertices) and the {\it bounded aspect ratio property} often called {\it regularity condition} (the diameter of the element and the diameter of the largest ball inscribed are comparable). 

Classical results on general Lagrange finite elements consider the regularity condition \cite{CR}. On triangular elements, the error estimate holds under the minimum angle condition \cite{Ze,Z}. However, on triangles, the minimum angle condition and the regularity condition are equivalent. From \cite{BA,BG,J} we know that the weakest sufficient condition on triangular elements is the maximum angle condition. Some examples can be constructed in order to show that if a family of triangles does not satisfy the maximum angle condition, then the error estimate on these elements does not hold. 

Recently, it was proved \cite{AM:2} that, for quadrilaterals elements, the minimum angle condition ($mac$) is the weakest known geometric condition required to obtain the classical $W^{1,p}$-error estimate, when $1 \leq p < 3$, to any arbitrary order $k$ greater than 1. Moreover, in this case, $mac$ is also necessary. In \cite{AM,AM:2} it was proved that the {\emph{double angle condition}} (any interior angle bounded away from zero and $\pi$) is a sufficient requirement to obtain the error estimate for any order and any $p \geq 1$.  When $k=1$ and $1 \leq p < 3$, a less restrictive condition ensures the error estimate \cite{AD,AM}: the {\it regular decomposition property} ($RDP$). Property $RDP$ requires that after dividing the quadrilateral into two triangles along one of its diagonals, each resultant triangle verifies the maximum angle condition and the quotient between the length of the diagonals is uniformly bounded. 

This brief picture intends to show that study of sharp geometric restrictions on finite elements under which the optimal error estimate remains valid is an interesting and active field of research.

In \cite{GRB,GRB:2}, geometric conditions on general polygons are given in order to guarantee the error estimate for interpolants built from generalized barycentric coordinates, and the question about identifying sharp geometric restrictions in this setting is proposed. In this work, we address the question for the first-order Wachspress interpolation operator. 

We show that the three sufficient conditions considered in \cite{GRB} ({\it regularity condition}, {\it maximum angle condition} and {\it minimum edge length property}) are actually equivalent to the last two since the regularity condition is a consequence of the maximum angle condition and the minimum edge length property. Then we exhibit families of polygons satisfying only one of these conditions and show that the interpolation error estimate does not hold to adequate functions. In this sense, the {\it maximum angle condition} and the {\it minimum edge length property} can be regarded as sharp geometric requirements to obtain the optimal error estimate. 

This work is structured as follows: In Section \ref{geoimpl}, we introduce notation and exhibit some basic relationships between different geometric conditions on general convex polygons. Section \ref{wach} is devoted to recall Wachspress coordinates and some elementary results associated to them; a brief picture about error estimates for the first-order Wachspress interpolation operator is also given there. Finally, in Section \ref{sharp}, we present two counterexamples to show that $MAC$ and $melp$ are sharp geometric requirements under which the optimal error estimate is valid.

\section{Geometric conditions}
\label{geoimpl}
\setcounter{equation}{0}

In order to introduce notation and formalize the requirements of each geometric condition, we give the following definitions. From now on, $\Omega$ will refer to a general convex polygon.

\medskip
\begin{enumerate}
	\item[(i)] {\it (Bounded aspect ratio property)} We say that $\Omega$ satisfies the {\emph{bounded aspect ratio property}} (also called {\emph{regularity condition}}) if there exists a constant $\sigma>0$ such that 
	\bq
	\label{barp}
	\frac{diam(\Omega)}{\rho(\Omega)} \le \sigma,
	\eq
	where $\rho(\Omega)$ is the diameter of the maximum ball inscribed in $\Omega$. In this case, we write $barp(\sigma)$.
	
	\item[(ii)] {\it (Minimum edge length property)} We say that $\Omega$ satisfies the {\emph{minimum edge length property}} if there exists a constant $d_m>0$ such that 
	\bq
	\label{mel}
	0<d_m \leq \frac{\left\| {\bv}_i-{\bv}_j \right\|}{diam(\Omega)}
	\eq
	for all $i \neq j$, where ${\bv}_1, {\bv}_2, \dots, {\bv}_n$ are the vertices of $\Omega$. In this case, we write $melp(d_m)$. 
	
	\item[(iii)] {\it (Maximum angle condition)} We say that $\Omega$ satisfies the {\emph{maximum angle condition}} if there exists a constant $\psi_M>0$ such that 
	\bq
	\label{MAC}
	\beta \leq \psi_M < \pi
	\eq
	for all interior angle $\beta$ of $\Omega$. In this case, we write $MAC(\psi_M)$,

	\item[(iv)] {\it (Minimum angle condition)} We say that $\Omega$ satisfies the {\emph{minimum angle condition}} if there exists a constant $\psi_m>0$ such that
	\bq
	\label{mac}
	0 < \psi_m \leq \beta.
	\eq
	for all interior angle $\beta$ of $\Omega$. In this case, we write $mac(\psi_m)$.
\end{enumerate}

All along this work, when we say {\it regular polygon}, we refer to a polygon satisfying the regularity condition given by \eqref{barp}.

\subsection{Some basic relationships}

It is well known that regularity assumption implies that the minimum interior angle is bounded away from zero. We state this result in the following lemma 

\begin{lemma}
	\label{lemma:regmac}
	If $\Omega$ is a convex polygon satisfying $barp(\sigma)$, then $\Omega$ verifies $mac(\psi_m)$ where $\psi_m$ is a constant depending only on $\sigma$.
\end{lemma} 

\proof See for instance \cite[Proposition 4 (i)]{GRB}. \qed

\medskip
Considering the rectangle $R=[0,1] \times [0,s]$, where $0<s<1$, and taking $s \to 0^+$, we see that the converse statement of Lemma \ref{lemma:regmac} does not hold. Indeed, $R$ verifies the $mac(\pi/2)$ (independently of $s$), but, when $s$ tends to zero, $R$ is not regular in the sense given by \eqref{barp}. However, on triangular elements, $barp$ and $mac$ are equivalent. We use this fact to show that, on general polygons, the regularity condition is a consequence of the minimum edge length property and the maximum angle condition. To our knowledge, this elementary result has not been established or demonstrated previously.

\begin{figure}[h]
	\resizebox{11.5cm}{5.2cm}{\includegraphics{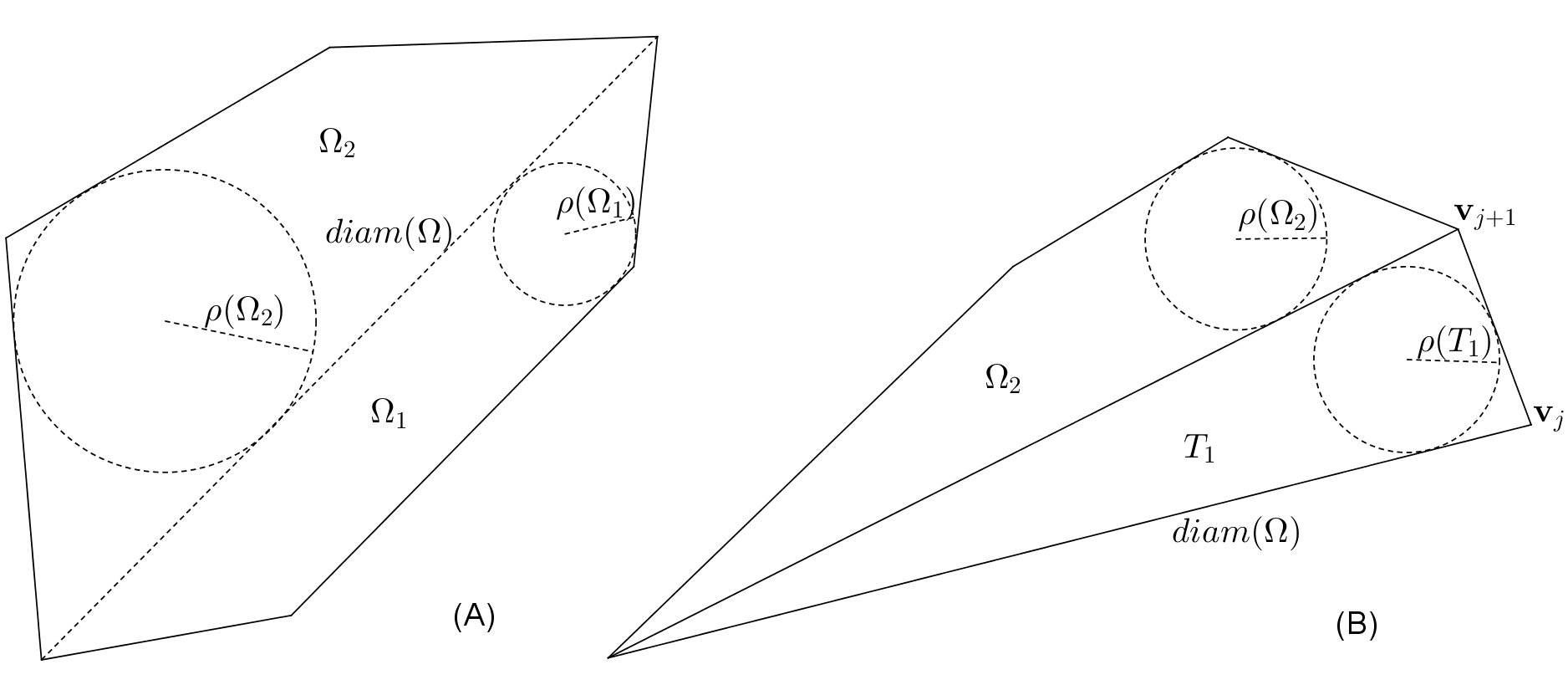}}	
	\caption{(A): A polygon with its diameter attained as the length of the straight line joining two non-consecutive vertices. (B): A polygon with its diameter attained as the length of the straight line joining two consecutive vertices.}
	\label{fig:macmel=reg}
\end{figure}

\begin{lemma}
	\label{lemma:macmelpeqreg}
	If $\Omega$ is a convex polygon satisfying $MAC(\psi_M)$ and $melp(d_m)$, then $\Omega$ verifies $barp(\sigma)$, where $\sigma=\sigma(\psi_M,d_m)$.
\end{lemma}

\proof We prove this by induction on the number $n$ of vertices of $\Omega$. If $n=3$, i.e., $\Omega$ is a triangle, the result follows from the law of sines. Indeed, we only have to prove that $\Omega$ has its minimum interior angle bounded away from zero. Let $\alpha$ be the minimum angle of $\Omega$ (if there is more than one choice, we choose it arbitrarily) and let $l$ be the length of its opposite side. Since $diam(\Omega)$ is attained on one side of $\Omega$, we can assume, without loss of generality, $l \neq diam(\Omega)$. We call $\beta$ the opposite angle to $diam(\Omega)$. It is clear that $\beta$ can not approach zero and since it is bounded above by $\psi_M$, we get that $1/\sin(\beta) \le C$ for some positive constant $C$. Then, from the law of sines and the assumption $melp(d_m)$, we have
$$\frac{\sin(\alpha)}{\sin(\beta)} = \frac{l}{diam(\Omega)} \geq d_m.$$

In consequence, $\sin(\alpha) \ge C^{-1} d_m$ which proves that $\alpha$ is bounded away from zero.

Let $n>3$. Since the diameter of $\Omega$ realizes as the length of its longest {\it diagonal}, i.e., the longest straight line joining two vertices of $\Omega$, we need to consider two cases depending if these vertices are consecutive or not. 

Assume that $diam(\Omega)$ is attained as the length of the line joining two non-consecutive vertices (these may not be unique, in this case we choose them arbitrarily). We can divide $\Omega$ by this diagonal into two convex polygons $\Omega_1$ and $\Omega_2$ with less number of vertices (see Figure \ref{fig:macmel=reg} (A)). It is clear that both of them satisfy $MAC(\psi_M)$ and, since $diam(\Omega_i)=diam(\Omega)$ and the set of vertices of $\Omega_i$ is a subset of the vertices of $\Omega$, we conclude that $\Omega_i$ also verifies $melp(d_m)$. Therefore, by the inductive hypothesis, $\Omega_1$ and $\Omega_2$ verify $barp(\sigma_1)$ and $barp(\sigma_2)$, respectively, for some constants $\sigma_1, \sigma_2$ depending only on $\psi_M$ and $d_m$. Then, since $\rho(\Omega) \geq \rho(\Omega_i)$, $i=1,2$, we have
$$\ds \frac{diam(\Omega)}{\rho(\Omega)} = \frac{diam(\Omega_i)}{\rho(\Omega)} \leq \frac{diam(\Omega_i)}{\rho(\Omega_i)} \leq \sigma_i.$$

Finally, if $diam(\Omega)$ is attained on a side of $\Omega$, i.e., is the length of the line joining two consecutive vertices ${\bv}_{j-1}$ and ${\bv}_j$ (these may not be unique, in this case we choose them arbitrarily), we divide $\Omega$ by the diagonal joining ${\bv}_{j-1}$ and ${\bv}_{j+1}$ into the triangle $T_1=\Delta({\bv}_{j-1}{\bv}_j{\bv}_{j+1})$ and a convex polygon $\Omega_2$ (see Figure \ref{fig:macmel=reg} (B)). It is clear that $T_1$ verifies $melp(d_m)$ and $MAC(\psi_M)$, so (by the case $n=3$) we have that $T_1$ satisfies $barp(\sigma_1)$ for some positive constant $\sigma_1$. Then, since $diam(T_1)=diam(\Omega)$ and $\rho(\Omega) \geq \rho(T_1)$, we have 
$$\ds \frac{diam(\Omega)}{\rho(\Omega)} = \frac{diam(T_1)}{\rho(\Omega)} \leq \frac{diam(T_1)}{\rho(T_1)} \leq \sigma_1.$$ \qed

\begin{cor}
	\label{cor:equiv}
	$[MAC, melp]$ and $[barp, MAC, melp]$ are equivalent conditions.
\end{cor}

Finally, notice that reciprocal statement of Lemma \ref{lemma:macmelpeqreg} is false. Consider the following families of quadrilaterals: $\mathcal{F}_1=\{ K(1,1-s,s,1-s) \}_{0<s<1}$ where $K(1,1-s,s,1-s)$ denotes the convex quadrilateral with vertices $(0,0), (1,0), (s, 1-s)$ and $(0,1-s)$, and $\mathcal{F}_2=\{ K(1,1,s,s) \}_{1/2<s<1}$ where $K(1,1,s,s)$ denotes the convex quadrilateral with vertices $(0,0), (1,0), (s, s)$ and $(0,1)$. Clearly, any quadrilateral belonging to $\mathcal{F}_1 \cup \mathcal{F}_2$ is regular in the sense given by \eqref{barp}. Each element of $\mathcal{F}_1$ satisfies $MAC(3\pi/4)$, but taking $s \to 0^+$, we see that the minimum edge length property is violated. On the other hand, each element of $\mathcal{F}_2$ verifies $melp(1/2)$; but taking $s \to 1/2^+$, we see that the maximum angle condition is not satisfied.

\section{Wachspress coordinates and the error estimate}
\label{wach}
\setcounter{equation}{0}

\subsection{Wachspress coordinates}

\begin{figure}[h]
	\centering
	\resizebox{11.6cm}{5cm}{\includegraphics{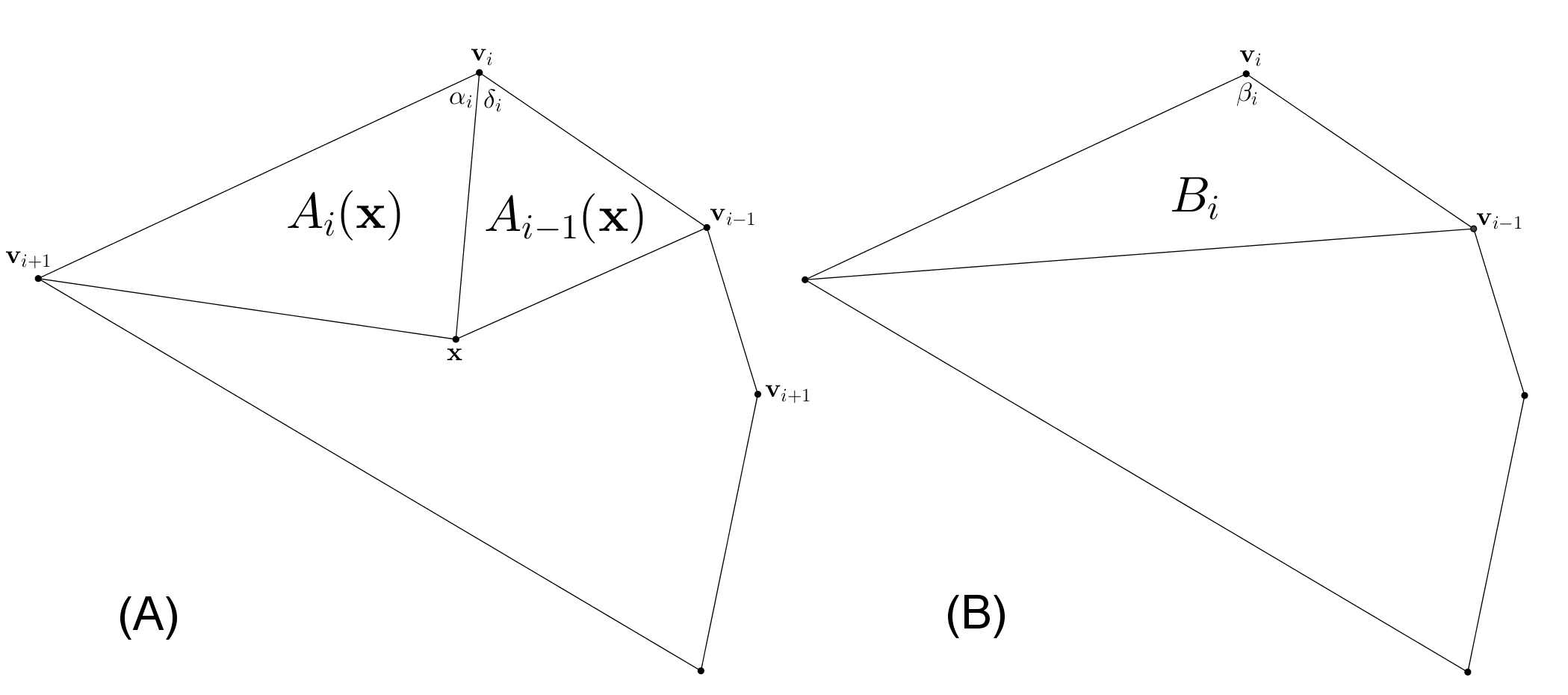}}	
	\caption{(A): Notation for $A_i({\bx})$. (B): Notation for $B_i$.}
	\label{fig:notation}
\end{figure}

We start this section by remembering the definition of Wachspress coordinates and some of their main properties \cite{Fl:2, W}. Henceforth, we denote by ${\bv}_1, {\bv}_2, \dots, {\bv}_n$ the vertices of $\Omega$ enumerated in counterclockwise order starting in an arbitrary vertex. Let $\bx$ denote an interior point of $\Omega$ and let $A_i(\bx)$ denote the area of the triangle with vertices $\bx$, ${\bv}_i$ and ${\bv}_{i+1}$, i.e., $A_i({\bx})=|\Delta ({\bx} {\bv}_i {\bv}_{i+1})|$, where, by convention, ${\bv}_0:= {\bv}_n$ and ${\bv}_{n+1}:={\bv}_1$. Let $B_i$ denote the area of the triangle with vertices ${\bv}_{i-1}$, ${\bv}_i$ and ${\bv}_{i+1}$, i.e., $B_i=|\Delta ({\bv}_{i-1} {\bv}_i {\bv}_{i+1})|$. We summarize the notation in Figure \ref{fig:notation}.

\medskip
Define the Wachspress weight function $w_i$ as the product of the area of the “boundary” triangle, formed by ${\bv}_i$ and its two adjacent vertices, and the areas of the $n-2$ interior triangles, formed by the point ${\bx}$ and the polygon's adjacent vertices (making sure to exclude the two interior triangles that contain the vertex ${\bv}_i$), i.e.,
\bq
\label{wi}
\ds w_i({\bx}) = B_i \prod_{j \neq i,i-1} A_j(\bx).
\eq
After applying the standard normalization, Wachspress coordinates are then given by
\bq
\label{lambdai}
\ds \lambda_i({\bx}) = \frac{w_i({\bx})}{\sum_{j=1}^n w_j({\bx})}.
\eq

An equivalent expression of \eqref{wi} for $w_i$ is given in \cite{Mey}; the main advantages of this alternative expression is that the result is easy to implement and it shows that only the edge $\overline{{\bx} {\bv}_i}$ and its two adjacent angles $\alpha_i$ and $\delta_i$ are needed (see Figure \ref{fig:notation} (A)). Indeed, $w_i$ can be written as 
\bq
\label{weights}
w_i({\bx}) = \frac{\cot(\alpha_i)+\cot(\delta_i)}{\left\| {\bx}-{\bv}_i \right\|^2}
\eq

where $\alpha_i=\angle\ {\bx} {\bv}_i {\bv}_{i+1}$ and $\delta_i=\beta_i-\alpha_i$ with $\beta_i$ being the inner angle of $\Omega$ associated to ${\bv}_i$ (see Figure \ref{fig:notation}). The evaluation of the Wachspress basis function is carried out using elementary vector calculus operations. The angles $\alpha_i$ and $\delta_i$ are not explicitly computed, as suggested in \cite{Mey}, vector cross product and vector dot product formulas are used to find the cotangents.

\medskip
Wachspress coordinates have the well-known following properties:
\bit
\item[(I)] {\it (Non-negativeness)} $\lambda_i \geq 0$ on $\Omega$.

\item[(II)] {\it (Linear Completeness)} for any linear function $\ell :\Omega \to \R$, there holds $\ell = \sum_{i} \ell({\bv}_i) \lambda_i$.

\item[] (Considering the linear map $\ell \equiv 1$ yields $\sum_{i} \lambda_i = 1$; this property is usually named {\it partition of unity}).

\item[(III)] {\it (Invariance)} If $L:\R^2 \to \R^2$ is a linear map and $S:\R^2 \to \R^2$ is a composition of rotation, translation and uniform scaling transformations, then $\lambda_i({\bx})=\lambda_i^L(L({\bx}))=\lambda_i^S(S({\bx}))$, where $\lambda_i^F(F({\bx}))$ denotes a set of barycentric coordinates on $F(\Omega)$.

\item[(IV)] {\it (Linear precision)} $\sum_{i} {\bv}_i \lambda_i({\bx})={\bx}$, i.e., every point on $\Omega$ can be written as a convex combination of the vertices ${\bv}_1, {\bv}_2, \dots, {\bv}_n$.

\item[(V)] {\it (Interpolation)} $\lambda_i({\bv}_j)=\delta_{ij}$.
\eit

\subsection{Error estimate to the first-order Wachspress interpolation operator} 

We only give a brief overview of some definitions and results which are of interest to us; for more details we refer to \cite{Das, GRB, Suk:2, Suk}.

Let $\{ \lambda_i \}$ be the Wachspress coordinates associated to $\Omega$ (see \eqref{lambdai}). Then, we can consider the first-order interpolation operator $I:H^2(\Omega) \to span \{ \lambda_i \} \subset H^1(\Omega)$ defined as
\bq
\label{defI}
\ds I_{\Omega}u=Iu := \sum_{i} u({\bv}_i) \lambda_i.
\eq

Properties (I)-(V) of the Wachspress coordinates (more generally, generalized barycentric coordinates) guarantee that $I$ has the desirable properties of an interpolant. For this interpolant, called here the {\it first-order Wachspress interpolation operator}, the optimal convergence estimate 
\bq
\label{errorestimate}
\left\| u-Iu \right\|_{H^1(\Omega)} \leq C diam(\Omega) |u|_{H^2(\Omega)}
\eq
on polygons satisfying $[barp, MAC, melp]$ was proved \cite[Lemma 6]{GRB}.

\begin{rem}
	\label{rem:red}
	Thanks to {\rm Corollary \ref{cor:equiv}}, we can affirm that \eqref{errorestimate} holds on general convex polygons satisfying $[MAC, melp]$.
\end{rem}

\section{About sharpness on geometric restrictions} 
\label{sharp}
\setcounter{equation}{0}

Since $[MAC, melp]$ are sufficient conditions to obtain \eqref{errorestimate}, we wonder if some of these requirements can be relaxed in order to obtain the error estimate. This question was partially answered in \cite{GRB}, where a counterexample, using pentagonal elements, is given in order to show that the $MAC$ can not be removed. For the sake of completeness, in Counterexample \ref{necmac}, we give a family of quadrilateral elements which does not satisfy $MAC$ but it verifies $melp$ and \eqref{errorestimate} does not hold. This example shows two things: $MAC$ is necessary in order to obtain the error estimate and, since every element in this family is regular in the sense given by \eqref{barp}, $barp$ is not enough to obtain \eqref{errorestimate}. 

On the other hand, in Counterexample \ref{necmel}, we present a family of quadrilaterals which does not satisfy $melp$ but it verifies $MAC$ and \eqref{errorestimate} does not hold. Then, in order to obtain the interpolation error estimate, $melp$ is necessary.

In this sense, the question raised in \cite{GRB} about identifying sharp geometric restrictions under which the error estimates for the first-order Wachspress interpolation operator holds can be considered as answered. 

\begin{figure}[h]
	\resizebox{12cm}{5cm}{\includegraphics{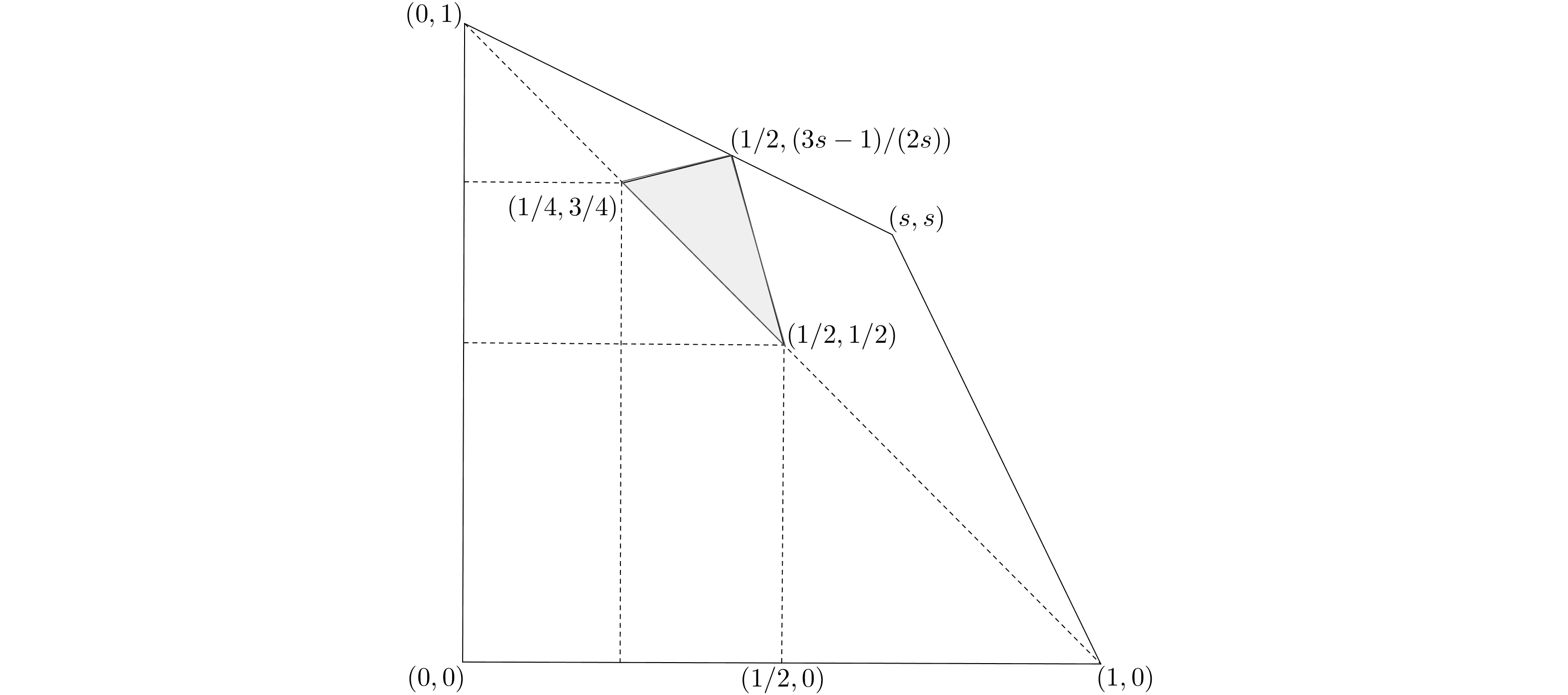}}
	\caption{Schematic picture of $K_s$ and $T_s$ (hatched area) considered in Counterexample \ref{necmac}.}
	\label{fig:cex1}
\end{figure}

\begin{cex}
	\label{necmac}
	Consider the convex quadrilateral $K_s$ with the vertices ${\bv}_1=(0,0), {\bv}_2=(1,0), {\bv}_3=(s,s)$ and ${\bv}_4=(0,1)$, where $1/2<s<1$. We will be interested in the case when $s$ tends to $1/2$ since then the family of quadrilaterals $\{ K_s \}$ does not satisfy the maximum angle condition although it satisfies $melp(1/2)$.
	
	Consider the function $u({\bx})=x(1-x)$. Since $u({\bv}_1)=0=u({\bv}_2)=u({\bv}_4)$, we have
	$$Iu({\bx})=u({\bv}_3) \lambda_3({\bx})= s(1-s) \lambda_3({\bx}).$$
	
	An straightforward computation yields
	$$\ds \lambda_3({\bx}) = \frac{(2s-1)x}{s} \frac{y}{(s-1)(x+y)+s},$$
	
	therefore
	$$\ds \frac{\partial \lambda_3}{\partial y} = \frac{(2s-1)x}{s} \frac{(s-1)x+s}{[(s-1)(x+y)+s]^2}.$$
	
	Consider the triangle $T_s$ with vertices $(1/4,3/4)$, $(1/2,1/2)$ and $(1/2, (3s-1)/(2s))$ {\rm (see Figure \ref{fig:cex1})}. Then, on $T_s$, we have $1/4 \leq x \leq 1/2$, $1/2 \le y \le (3s-1)/(2s)$ and $x+y \geq 1$, so it follows that
	$$0<(s-1)(x+y)+s \leq 2s-1 
	\quad \text{and} \quad
	(s-1)x+s \geq (3s-1)/2$$
	
	and hence
	$$\ds \frac{\partial \lambda_3}{\partial y} \geq \frac{(2s-1)}{4s} \frac{3s-1}{2(2s-1)^2}=\frac{3s-1}{8s(2s-1)}.$$
	
	Then
	$$|u-Iu|_{H^1(K_s)} \ge \left\| \frac{\partial (u-Iu)}{\partial y} \right\|_{L^2(K_s)} = \left\| \frac{\partial Iu}{\partial y} \right\|_{L^2(K_s)}	= s(1-s)\left\| \frac{\partial \lambda_3}{\partial y} \right\|_{L^2(K_s)}$$
	and, consequently,
	$$|u-Iu|_{H^1(K_s)} \ge s(1-s)\left\| \frac{\partial \lambda_3}{\partial y} \right\|_{L^2(T_s)}.$$
	Since $|T_s|=(2s-1)/(2^4s)$, we have
	$$\left\| \frac{\partial \lambda_3}{\partial y} \right\|_{L^2(T_s)}^2 \geq \frac{(3s-1)^2}{(8s)^2(2s-1)^2}|T_s|=
	\frac{(3s-1)^2}{2^{10}s^3(2s-1)} \to \infty$$
	when $s \to 1/2^+$. Finally, as $|u|_{H^2(K_s)} = 2 |K_s|^{1/2} \leq 2$ and $diam(K_s)=\sqrt{2}$, we conclude that \eqref{errorestimate} can not hold.
\end{cex}

\begin{figure}[h]
	\centering
	\resizebox{11.3cm}{5cm}{\includegraphics{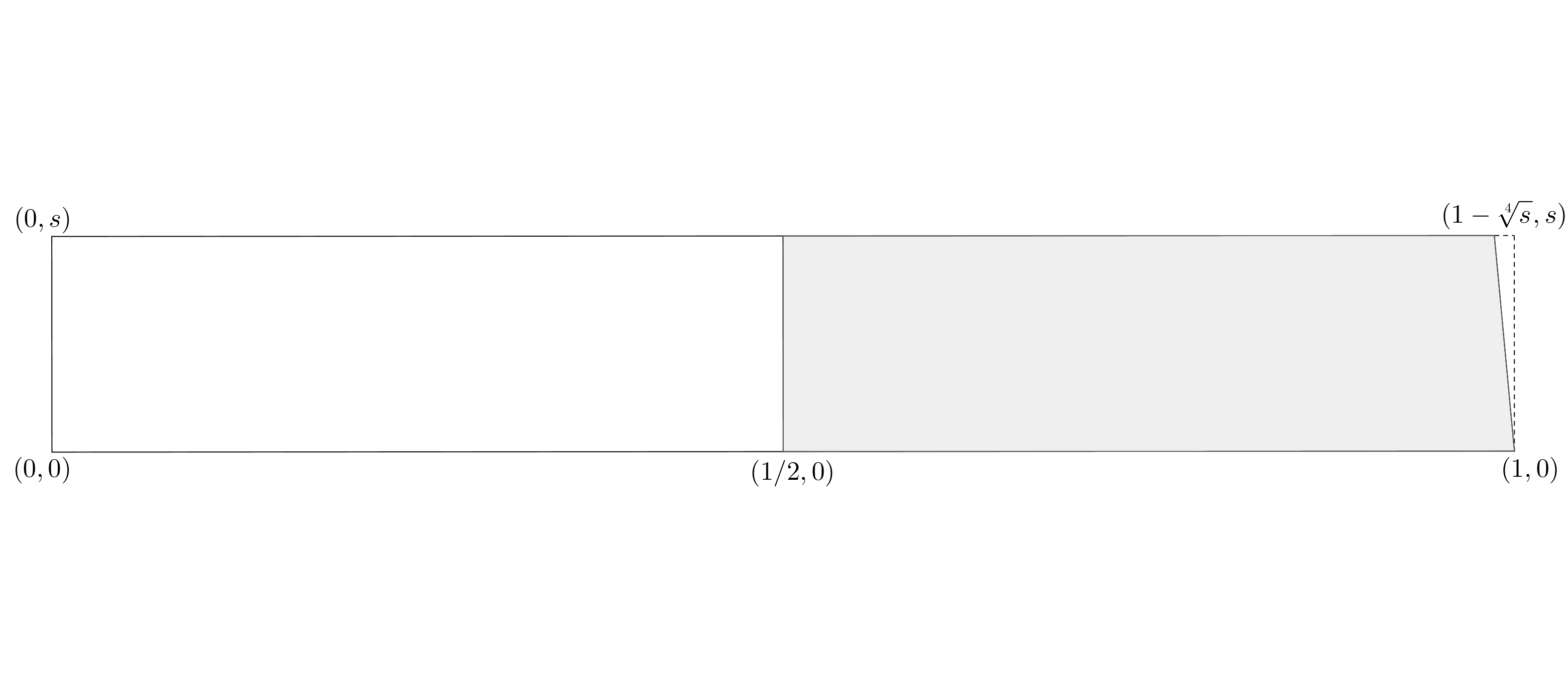}}
	\caption{Schematic picture of $K_s$ and $D_s$ (hatched area) considered in Counterexample \ref{necmel}.}
	\label{fig:cex2}
\end{figure}

\begin{cex}
	\label{necmel}
	Consider now the convex quadrilateral $K_s$ with the vertices ${\bv}_1=(0,0), {\bv}_2=(1,0), {\bv}_3=(1-\sqrt[4]{s},s)$ and ${\bv}_4=(0,s)$, where $0 < s < (1/2)^4$. Note that the family of quadrilaterals $\{ K_s \}$ satisfies $MAC(\pi/2+\tan^{-1}(2^3))$ $($independently of $s)$ but it does not satisfy the minimum edge length property when $s$ tends to zero since $\left\| {\bv}_1-{\bv}_4 \right\| = s \to 0^+$ and $diam(K_s) \sim 1$.
	
	Consider the function $u({\bx})=x^2$. Since $u({\bv}_1)=0=u({\bv}_4)$, we have, calling $a := 1-\sqrt[4]{s}$,
	$$Iu({\bx})=u({\bv}_2) \lambda_2({\bx})+u({\bv}_3) \lambda_3({\bx}) = \lambda_2({\bx})+ a^2 \lambda_3({\bx})$$
	
	where
	$$\lambda_2({\bx})=\frac{x(s-y)}{s+y(a-1)}  \quad \text{and} \quad 
	\lambda_3({\bx})=\frac{xy}{s+y(a-1)}.$$
	
	A simple computation yields
	$$\frac{\partial (Iu-u)}{\partial y} = \frac{\partial Iu}{\partial y} = \frac{xsa(a-1)}{(s+y(a-1))^2}.$$
	
	Let $D_s = K_s \cap \{ x \geq 1/2 \}$ $($see {\rm Figure \ref{fig:cex2})}. Since $a-1 <0$, we get $s+y(a-1) \leq s$ and then, on $D_s$, we have 
	$$\left| \frac{\partial (Iu-u)}{\partial y} \right| \geq \frac{xa(1-a)}{s} \geq \frac{a(1-a)}{2s}.$$ 
	
	Therefore, 
	$$|Iu-u|_{H^1(K_s)}^2 \geq
	\left\| \frac{\partial (Iu-u)}{\partial y} \right\|_{L^2(K_s)}^2 \geq 
	\left\| \frac{\partial (Iu-u)}{\partial y} \right\|_{L^2(D_s)}^2 \geq
	\frac{a^2(1-a)^2}{4s^2} |D_s|,$$
	
	and since $|D_s|=as/2$, we conclude that
	$$|Iu-u|_{H^1(K_s)}^2 \geq 
	\frac{a^3(1-a)^2}{8s} =
	\frac{(1-\sqrt[4]{s})^3}{8\sqrt{s}}$$
	
	which tends to infinity when $s$ tends to zero. Finally, since $|u|_{H^2(K_s)} = 2 |K_s|^{1/2} \leq 2$ and $diam(K_s) \sim 1$, we conclude that \eqref{errorestimate} can not hold.
\end{cex}


\begin{thebibliography}{153}
	
	\bibitem{AD} Acosta G., Dur\'an R. G. (2000):
	Error estimates for $Q_1$ isoparametric elements satisfying a weak
	angle condition, SIAM J. Numer. Anal., 38, 1073-1088.
	
	\bibitem{AM} Acosta G., Monz\'on G. (2006):
	Interpolation error estimates in $W^{1,p}$ for degenerate $Q_1$
	isoparametric elements, Numer. Math., 104, 129-150.
	
	\bibitem{AM:2} Acosta G., Monz\'on G. (2017):
	The minimal angle condition for quadrilateral  finite elements
	of arbitrary order. J. Comput. Appl. Math., 317, 218-234.
	
	\bibitem{BA} Babu$\breve{\mbox{s}}$ka I., Aziz A. K. (1976): 
	On the angle condition in the finite element method. 
	SIAM J. Numer. Anal., 13, 214-226.
	
	\bibitem{BG} Barnhill R. E., GregoryJ. A. (1976): 
	Sard kernel theorems on triangular domains with applications to finite element error bounds. Numer. Math., 25, 215-229.
	
	\bibitem{CR} Ciarlet P. G., Raviart P. A. (1978):
	The finite element method for elliptic problems. 
	North-Holland Publishing Company. Amsterdam.
	
	\bibitem{Das} Dasgupta G. (2003):
	Interpolants within Convex Polygons: Wachspress' Shape Functions. 
	Journal of Aerospace Engineering, 16, 1-8.
	
	\bibitem{Fl:2} Floater, M., Hormann, K., K\'os, G. (2006): 
	A general construction of barycentric coordinates over
	convex polygons. Advances in Computational Mathematics, 24, 311-331. 
	
	\bibitem{GRB} Gillette, A., Rand, A., Bajaj, C. (2012): 
	Error estimates for generalized barycentric interpolation. 
	Adv. Comput. Math., 37, 417-439. 
	
	\bibitem{GRB:2} Gillette, A., Rand, A., Bajaj, C. (2013): 
	Interpolation error estimates for mean value coordinates over convex polygons. 
	Adv. Comput. Math., 39, 327-347. 
	
	\bibitem{J} Jamet P. (1976): 
	Estimation de l'erreur pour des \'el\'ements finis droits presque d\'eg\'en\'er\'es.
	RAIRO Anal. Num\'er., 10, 43-60.
	
	\bibitem{Mey} Meyer M., Lee H., Barr A. , Desbrun M. (2002):
	Generalized barycentric coordinates on irregular polygons. Journal of Graphics Tools, 7, 13-22.
	
	\bibitem{Suk} Sukumar N., Tabarrai A. (2004): 
	Conforming polygonal finite elements. Int J. Numer. Methods Eng., 61, 2045-2066.
	
	\bibitem{Suk:2} Sukumar N., Malsch E. A. (2006):
	Recent Advances in the Construction of Polygonal Finite Element Interpolants. 
	Arch. Comput. Meth. Engng., 13, 129-163.
	
	\bibitem{W} Wachspress E. L. (1975): 
	A rational finite elements basis. Academic press, New York. 
	
	\bibitem{Ze} Zenisek A. (1969): 
	Convergence of the finite element method for boundary value problems of a system of elliptic equations. 
	Apl. Mat., 14, 355-377.
	
	\bibitem{Z} Zl\'amal M (1968): 
	On the finite element method. 
	Numer. Math., 12, 394-409.
	
\end{thebibliography}
\end{document}